
\input amstex
\magnification=\magstep 1

\TagsOnRight
\documentstyle{amsppt}
\overfullrule=0pt \vsize=550pt

\define\bsk{\bigskip}
\define\nli{\newline}

\define\DSM{Dobrushin and Shlosman's strong mixing conditions}
\define\Duc{Dobrushin's uniqueness condition}
\define\SMW{strong mixing condition}
\define\SM{strong mixing condition }
\define\cond{condition}

\define\wre{with respect to}
\define\sqe{sequence}
\define\sat{satisfying}
\define\sats{satisfies}
\define\sy{satisfy}

\define\ms{measure}

\define\ineq{inequality}
\define\ineqs{inequalities}

\define\LSI{logarithmic Sobolev \ineq}
\define\LSIs{logarithmic Sobolev inequalities}
\define\LSc{logarithmic Sobolev constant}
\define\relent{relative entropy}

\define\Mk{Markov kernel}
\define\pms{probability measure}

\define\rv{random variable}
\define\dt{distribution}

\define\cd{conditional distribution}

\define\jd{joint distribution}

\define\lsp{local specification}
\define\config{configuration}

\define\pa{\partial}

\define\dis{distance}

\define\coord{coordinate}

\define\X{\Cal X}
\define\Lw{\Cal L}

\define\I{\Cal I}

\define\Zd{\Bbb Z^d}
\define\Id{\Bbb I\Bbb d}

\define\PP{\Cal P}
\define\PXn{\PP(\X^n)}
\define\PXLa{\PP(\X^\La)}
\define\byL{\bar y_\La}

\define\a{\alpha}

 \define\la{\lambda}

\define\f{\varphi}

\define\Ga{\Gamma}

\define\de{\delta}
\define\ga{\gamma}
\define\La{\Lambda}

\define\sumi{\sum_{i=1}^n}

\define\qute{\quad\text}
\define\qand{\quad\text{and}\quad\ }

\define\setm{\setminus}
\define\su{\subset}
\define\susu{\subset\subset}

\document

\topmatter

\title
{Logarithmic Sobolev inequalities in  discrete product spaces:
a proof by a
 transportation cost distance}
\endtitle

\leftheadtext {Logarithmic Sobolev inequality }
\rightheadtext {Proof by a transportation cost distance}

\author Katalin Marton \endauthor
\affil MTA Alfr\'ed R\'enyi  Institute of Mathematics
\endaffil

\date{July 3, 2015}
\enddate

\address H-1364 POB 127, Budapest, tel. 36 (1) 4838300, Hungary
\endaddress
\email marton\@renyi.hu
\endemail
\footnote{This work was supported in part by the grants OTKA  K 105840
of the
Hungarian Academy of Sciences}

\abstract
The aim of this paper is to prove
an \ineq\ between relative entropy and the sum of average conditional
relative entropies of the following form:
For a fixed \pms\ $q^n$ on $\X^n$, ($\X$ is a finite set), and
any \pms\ $p^n=\Lw(Y^n)$ on $\X^n$
$$
\align
&D(p^n||q^n)\leq\\
& Const.
\sum_{i=1}^n
\Bbb E_{p^n}
D(p_i(\cdot|Y_1,\dots, Y_{i-1},Y_{i+1},\dots, Y_n)
||
  q_i(\cdot|Y_1,\dots, Y_{i-1},Y_{i+1},\dots, Y_n)),
  \tag *
\endalign
$$
where $p_i(\cdot|y_1,\dots, y_{i-1},y_{i+1},\dots, y_n)$ and
$q_i(\cdot|x_1,\dots, x_{i-1},x_{i+1},\dots, x_n)$ denote
the local specifications for $p^n$ resp. $q^n$, i.e.,
the \cd s of the $i$'th
coordinate, given the other coordinates.
The constant shall depend on the properties of the local specifications of $q^n$.
\bsk

The inequality (*)  is meaningful  in product spaces,
both in the discrete and the continuous case, and
 can be used to prove a
\LSI\ for $q^n$, provided uniform \LSIs\ are available for
\newline
$q_i(\cdot|x_1,\dots, x_{i-1},x_{i+1},\dots, x_n)$, for all fixed $i$
and  all fixed $(x_1,\dots, x_{i-1},x_{i+1},\dots, x_n)$.
 (*) directly
implies that
 the Gibbs sampler
associated with $q^n$ is a contraction for relative entropy.

\bsk

In this paper we derive \ineq\ (*), and thereby a  \LSI, in
discrete
product spaces,
by proving \ineqs\ for an appropriate Wasserstein-like \dis.
\bsk

A \LSI\ is, roughly speaking, a contractivity property of relative entropy
\wre\ some Markov semigroup.
It is much easier to prove contractivity for a \dis\
  between \ms s, than for \relent,
  since \dis s \sy\ the triangle \ineq, and
  for them
well known linear tools, like estimates through matrix norms
can be applied.
\bsk

\endabstract

\endtopmatter

\beginsection 1. Introduction  and statement of some results.

Let $\X$ be a finite set,
and $\X^n$ the set of $n$-length \sqe s from  $\X$. Denote by $\PXn$
the space of  \pms s on $\X^n$.
For a \sqe\ $ x^n\in \X^n$ we denote by $x_i$ the $i$-th \coord\ of  $ x^n$.

\bsk

We consider a reference \pms\ $q^n\in \PXn$ which will be fixed throughout
Sections 1-3. In  section 4 we still  consider a fixed \pms\
denoted by $q$, with some subscript.
\bsk

The aim of this paper is to prove \LSIs\ for  \ms s
on discrete
product spaces,
by proving \ineqs\ for an appropriate Wasserstein-like \dis.
A \LSI\ is, roughly speaking, a contractivity property of relative entropy
\wre\ some    Markov semigroup.
It is much easier to prove contractivity for a \dis\
  between \ms s, than for \relent, since a \dis\ is symmetric and
  \sats\ the triangle \ineq. Our method shall be used to prove \LSIs\
for \ms s \sat\ a version
of Dobrushin's uniqueness condition, as well as Gibbs \ms s \sat\ a   \SM.
\bsk

To explain the results, we need
some definitions and some notation.
\bsk

\definition{Notation}
If $r$ and $s$ are two probability \ms s (on any measurable space)
then we denote by $|r-s|$ their variational distance:
$$
|r-s|=\sup_A\bigl|r(A)-s(A)\bigr|.
$$
\enddefinition
\bsk

\definition{Definition: $W_2$ distance} (c.f. [B-L-M],  Theorem 8.2)
\newline
For \pms s   $r^n,s^n\in \PXn$
let  $Z^n$ and $U^n$ represent  $r^n$ resp. $s^n$, i.e.,
 $Z^n$ and $U^n$ are \rv s with \dt s  $\Lw(Z^n)=r^n$ and $\Lw(U^n)=s^n$, respectively.
We    define
$$
W_2(r^n,s^n)=\min_\pi
\sqrt{\sumi Pr_\pi^2\bigl\{Z_i\neq U_i\bigr\}},
$$
where the minimum is taken over all \jd s $\pi=\Lw(Z^n,U^n)$ with marginals
$r^n$ and $s^n$.
\enddefinition
\bsk

Note that $W_2$ is a distance on $\PXn$,
but it cannot be
defined by taking the minimum expectation of a distance (or some power of a distance)
on  $\X^n$.
\bsk

\definition{Definition: Relative entropy, conditional relative entropy}
For \pms s $r$ and $s$ defined on a finite set $\Cal Z$,
we denote by
$D(r\Vert s)$ the relative entropy  of  $r$ \wre\   $s$:
$$
D(r\Vert s) = \sum_{u\in\Cal Z} r(u)\log \frac {r(u)} {s(u)},
$$
with the convention $0\log 0=0$ and $a\log 0=\infty$ for $a>0$.
If $Z$ and $U$ are \rv s with values in $\Cal Z$ and
distributed according to $r=\Lw(Z)$ resp. $s=\Lw(U)$,
then we shall also use the notation
$D(Z\Vert U)$ for the  relative entropy $D(r\Vert s)$.
If, moreover, we are
given a \pms\  $\pi=\Lw(S)$  on another finite set $\Cal S$, and
\cd s $\mu(\cdot|s)=\Lw(Z|S=s)$, $\nu(\cdot|s)=\Lw(U|S=s)$ then
we consider the average relative entropy
$$
\Bbb E_{\pi}D\bigl(\mu(\cdot|S)\Vert \nu(\cdot|S)\bigr)
=
\sum_{s\in \Cal S}
\pi(s)
D\bigl(\mu(\cdot|s)\Vert \nu(\cdot|s)\bigr).
$$

For $\Bbb E_{\pi}D(\mu(\cdot|S)\Vert \nu(\cdot|S))$
we shall   use either of the notations
$$
 D\bigl(\mu(\cdot|S)\Vert \nu(\cdot|S)\bigr),
\quad D\bigl(\mu(\cdot|S)\Vert U|S\bigr),
\quad D\bigl(Z|S)\Vert \nu(\cdot|S)\bigr),
\quad D\bigl(Z|S)\Vert U|S)\bigr)
$$
(omitting the symbol of expectation as is usual in information theory).
\enddefinition
\bsk

\definition{Notation}
\nli
For $y^n=(y_1,y_2,\dots,y_n)\in \X^n$ and $I\su[1,n]$,
we write
$$
y_I=(y_k: k\in I)\qand \bar y_I=(y_k: k\notin I).
$$
Moreover, if
 $p^n=\Lw(Y^n)$ then
$$
p_I\triangleq  \Lw(Y_I),\quad
p_I\bigl(\cdot|\bar y_I\bigr)\triangleq  \Lw\bigl(Y_I\bigm|\bar Y_I=\bar y_I\bigr),
\quad\bar p_I\triangleq  \Lw(\bar Y_I),
\quad\bar p_I(\cdot|y_I)\triangleq  \Lw(\bar Y_I|Y_I=y_I).
$$
If $I=\{i\}$ then we write $i$ instead of $\{i\}$.
\enddefinition
\bsk

\definition{Definition}
The \cd s $q_i(\cdot|\bar x_i)$  are called the \lsp s of
the \dt\ $q^n$.
\enddefinition
\bsk

\proclaim{Theorem 1}
\nli
Set
$$
\a\triangleq\min q_i(x_i|\bar x_i),
\tag 1.1
$$
where the minimum is taken over all  $x^n\in \X^n$ \sat\  $q(x^n)>0$ and all
$i\in [1,n]$.
Fix a $p^n=\Lw(Y^n)\in \PXn$ \sat
$$
q^n(x^n)=0 \implies p^n(x^n)=0.
\tag 1.2
$$
Assume that
 $q^n\in \PXn$
\sats\ all
the \ineqs\
$$
W_2^2\bigl(  p_I\bigl(\cdot|\bar y_I\bigr),
q_I\bigl(\cdot|\bar y_I\bigr)\bigr)
\leq
C\cdot
\Bbb E
\biggl\{
\sum_{i\in I}
\bigl|p_i(\cdot|\bar Y_i)-q_i(\cdot|\bar Y_i)\bigr|^2
\biggm|\bar Y_I=\bar y_I
\biggr\},
\tag 1.3
$$
where   $I\su [1,n]$ and  $\bar y_I\in \X^{[1,n]\setm I}$ is a fixed
\sqe. Then
$$
\align
&D(p^n||q^n)
\leq
\frac {4C} \a
\cdot
\sumi \Bbb E \bigl|p_i(\cdot|\bar Y_i)-Q_i(\cdot|\bar Y_i)\bigr|^2\\
&\leq
\frac {2C} \a
\cdot
\sumi D\bigl(Y_i|\bar Y_i||Q_i(\cdot|\bar Y_i)\bigr).
\tag 1.4
\endalign
$$
\endproclaim
\bsk

(Condition (1.2) is necessary, since otherwise $D(p^n||q^n)$ could be $\infty$, while
the middle term is always finite. On the other hand, for the \ineq\
between the  first and last terms it is not necessary
to assume (1.2), since if $D(p^n||q^n)=\infty$ then
the last term is $\infty$ as well.)
\bsk

\remark{Remark}
In [M] a  bound, analogous to the one relating
the first and  last terms of (1.4),
was proved  for  \ms s on Euclidean spaces. (Under reasonable conditions.)
 That bound   was used to derive
a \LSI, improving on an earlier result in [O-R].
In the present paper a \LSI\ shall be deduced from the first
\ineq\ in (1.4) (Corollary 2 to  Theorem 1).
\endremark
\bsk

Theorem 1 implies that the Gibbs sampler
(or Glauber dynamics) defined by the \lsp s of $q^n$
is a strict contraction for relative entropy.
\bsk

\definition{Definition: Gibbs sampler}
\nli
For $i\in[1,n]$ let
$\Ga_i: \PXn\mapsto  \PXn$ be the  \Mk\
$$
\Ga_i(z^n|y^n)
=
\de(\bar y_i,\bar z_i)
\cdot
q_i(z_i|\bar y_i ),\qquad y^n,z^n\in\X^n.
$$
(I.e., $\Ga_i$ leaves all, but the $i$-th, \coord s unchanged, and updates
the $i$-th \coord\ according to $q_i(y_i|\bar y_i )$.)
Finally, set
$$
\Ga
=\frac 1 n \cdot \sumi \Ga_i.
$$
I.e.,  $\Ga$  selects an $i\in[1,n]$ at random, and applies
$\Ga_i$.
It is easy to see that $\Ga$ preserves, and
is reversible \wre, $q^n$.
 $\Ga$ is called the Gibbs sampler governed by the \lsp  s of $q^n$.
\enddefinition

\bsk

\proclaim{Corollary 1 to  Theorem 1}
\newline
If $q^n$ on $\X^n$ satisfies  the conditions of Theorem 1
 then
$$
D(p^n\Gamma||q^n)
\leq
\biggl(1-\frac {\a}  {2nC} \biggr)
\cdot D(p^n||q^n).
\tag 1.5
$$
\endproclaim
\bsk

(1.5) follows from  Theorem 1 by the
\ineq\
$$
D(p^n\Gamma||q^n)\leq \frac 1 n\sumi  D(p^n\Gamma_i||q^n)
$$
(a consequence of the convexity of relative entropy), together with
 the identity
$$
D(p^n||q^n)- D(p^n\Gamma_i||q^n)
=
D\bigl(p_i\bigl(\cdot|\bar Y_i\bigr)
       ||q_i(\cdot|\bar Y_i)\bigr).
$$
\bsk

Theorem 1 also implies  Gross' \LSI\ defined as follows:
\bsk

\definition{Definition: \LSI\ for a \Mk}
\nli
Let $(\Cal Z, \pi)$ be a finite probability space, and
 $G: \Cal Z\mapsto\Cal Z$  a \Mk\ with invariant \ms\ $\pi$.
 The Dirichlet form associated with $G$ is
$$
\Cal E_G(f,f)=\bigl\langle\bigl(\Id-G\bigr)f,  f\bigr\rangle_\pi.
$$
We say that $G$ \sats\ a \LSI\ with \LSc\ $c$
if: for every
\pms\ $p$ on $\Cal Z$ we have
$$
D(p||\pi)
\leq
c\cdot \Cal E_G\bigl(\sqrt{f},\sqrt{f}\bigr),
$$
where $f(z)=p(z)/\pi(z)$.
\enddefinition
\bsk

The property expressed by the \LSI\ was  defined by L. Gross [Gr] in 1975.
For an introduction to \LSIs\ and their manifold interpretations and uses,
c.f. [L] and [R].
\bsk

Theorem 1 implies Gross' \LSI\ for the Gibbs sampler $\Ga$.
A simple calculation shows that
$$
\Cal E_\Ga
\biggl(
\sqrt{\frac {p^n} {q^n}},\sqrt{\frac {p^n} {q^n}}
\biggr)
=
\frac 1 n\cdot
 \Bbb E\sumi
\biggl(
1-\biggl(\sum_{y_i\in \X}\sqrt{p_i\bigl(y_i|\bar Y_i\bigr)
\cdot q_i\bigl(y_i|\bar Y_i\bigr)}\biggr)^2
\biggr).
$$
(Using the fact that, for fixed $\bar y_i$, the \ms\ $p^n\Ga_i$ does not depend on $y_i$, we just
 calculate the Dirichlet form for a matrix with identical rows.)
\bsk

\proclaim{Corollary 2 to  Theorem 1}
\newline
If $q^n$ on $\X^n$ satisfies  the conditions of Theorem 1
 then
$$
\align
&\frac 1 n\cdot D(p^n\Gamma||q^n)
\leq
\frac {4C} {\a n}
\cdot
\sumi \Bbb E
\biggl(
1-\biggl(\sum_{y_i\in \X}\sqrt{p_i\bigl(y_i|\bar Y_i\bigr)
\cdot q_i\bigl(y_i|\bar Y_i\bigr)}\biggr)^2
\biggr)\\
&=
\frac {4C} \a
\cdot
\Cal E_\Ga
\biggl(
\sqrt{\frac {p^n} {q^n}},\sqrt{\frac {p^n} {q^n}}
\biggr).
\endalign
$$
This can be considered a dimension free \LSI, since $\Ga$
only updates one coordinate.
\endproclaim
\bsk

Corollary 2 follows from  Theorem 1 by the following
\bsk

\proclaim{Lemma 1} (The proof is in  Appendix A)
\nli
Let $r$ and $s$ be two \pms s on $\X$.
Then
$$
|r-s|^2
\leq
1-\bigl(\sum_{y\in \X}\sqrt{r(y)s(y)}\bigr)^2.
$$
\endproclaim
\bsk

Theorem 1 can be applied to \dt s $q^n$ \sat\ the following  version
of Dobrushin's uniqueness condition:
\bsk

\definition{Definition: Dobrushin's uniqueness condition}
\newline
We say that $q^n$ satisfies (an $\Bbb L_2$-version of) Dobrushin's uniqueness condition with coupling matrix
$$
A=\bigl(a_{k,i}\bigr)_{k,i=1}^n,
$$
if: for any   pair of integers $k,i\in[1,n], k\neq i$ and any
two \sqe s $z^n,s^n\in \X^n$, differing only in the $k$'th coordinate, $$
|q_i(\cdot|\bar z_i)-q_i(\cdot|\bar s_i)|
\leq
a_{k,i},
\tag 1.6
$$
and, setting $a_{i,i}=0$ for all $i$,
$$
||A||_2<1.
$$
\enddefinition
\bsk

This differs from Dobrushin's original uniqueness condition where the norm   $||A||_1$ is assumed to be $<1$.
\bsk

\proclaim{Theorem 2}
\nli
Assume that the measure $q^n$ on $\X^n$ satisfies  Dobrushin's uniqueness condition with coupling matrix
$A$, $||A||_2<1$.
Then the conditions of Theorem 1 are satisfied with
$C= 1/ {\bigl(1-||A||\bigr)^2}$.
Thus for any $p^n\in \PXn$, \sat\ (1.3):
$$
\align
&D(p^n||q^n)
\leq
\frac 4 \a \cdot
\frac 1 {\bigl( 1-||A||\bigr)^2}
\cdot
\sumi \Bbb E \bigl|p_i(\cdot|\bar Y_i)-Q_i(\cdot|\bar Y_i)\bigr|^2\\
&\leq
\frac 2  \a \cdot
\frac 1 {\bigl( 1-||A||\bigr)^2}
\cdot
\sumi D\bigl(Y_i|\bar Y_i||Q_i(\cdot|\bar Y_i)\bigr),
\tag 1.7
\endalign
$$
and
$$
D(p^n\Gamma||q^n)
\leq
\biggl(1-\frac 1 n
\cdot
\frac \a  2
\cdot
\bigl( 1-||A||\bigr)^2    \biggr)
\cdot D(p^n||q^n).
\tag 1.8
$$
\endproclaim
\bsk

\remark{Remark}
In [Z] a \LSI\ is proved for discrete spin systems, where the title
suggests that it uses \Duc. However, the condition used there is reminiscent but not identical to \Duc. Moreover, an
\ineq\ of the form relating
the first and  last terms of (1.4) has been recently proved
in [C-M-T], assuming  conditions slightly reminiscent of \Duc.
\endremark
\bsk

Theorem 1 is proved in Section 2, and Theorem 2 in Section 3.
\bsk

In Section 4 we are going to deduce a \LSI\ from a \SMW,
 for \ms s $q$ on $\X^{\Zd}$.
 (Under the additional \cond\ that the \lsp s
$q_k(x_k|x_i, i\neq k)$, if not equal to $0$, are bounded from below.)
The \SMW\  we use is the same as \DSM, but we do not assume that $q$
is a Markov field. Our \SMW\  can also be considered as a generalization
of $\Phi$-mixing for (stationary) \pms s on $\X^{\Bbb Z}$.
 For  non-Markov stationary \pms s on $\X^{\Bbb Z}$ it is more restrictive
than usual strong mixing.
\bsk

The first  proof for the implication that  \DSM\ imply a \LSI\  for Markov fields
was given by  D. Stroock and B. Zegarlinski  [S-Z1],  [S-Z2] in 1992
(where the authors also proved the converse implication, i.e.
that \DSM\ for Markov fields are equivalent to the \LSI).
The arguments in
[S-Z2]  are quite hard to follow. In 2001, F. Cesi
proved that  \DSM\  imply a \LSI; his approach
is quite different from the previous ones, and much simpler.

\bsk

We feel that there is still room for alternative and perhaps simpler
  proofs in this important topic.
Moreover, our proof is valid without the Markovity assumption.
(It may be, though, that the proofs in [S-Z2] and [C] can also be
 generalized for the non-Markovian case, just it has not been tried.)
\bsk

We believe
 that the separate
parts of our proof  (Theorem 1 and the applicability of  Theorem 1)
are  comprehensible in themselves, thus
making the whole proof easier to follow.
\bsk
\bsk

\beginsection 2. Proof of Theorem 1.

We need the following

\proclaim{Lemma 2}
\newline
Let $r$ and $s$ be two \pms s on $\X$. Set
$$
\a_s=\min_{s(x)\neq 0} s(x).
$$
If $D(r||s)<\infty$ then
$$
D(r||s)\leq \frac 4 {\a_s} \cdot |r-s|^2.
\tag 2.1
$$
\endproclaim
\bsk

\remark{Remark}
\nli
Inequality (2.1) can be considered as a converse
to the Pinsker-Csisz\'ar-Kullback \ineq\ which says that
$$
|r-s|^2\leq \frac 1 2 D(r||s).
$$
However, there is no uniform converse: the reverse
\ineq\  must depend on $s$.
\endremark
\bsk

\demo{Proof}
\newline
Set $\X_+=\{x\in\X: s(x)>0\}$.  The following \ineq\ is well known:
$$
D(r||s)\leq \sum_{\X_+}\frac {|r(x)-s(x)|^2} {s(x)}.
$$
It follows that
$$
D(r||s)\leq
\frac 1 {\a_s}\cdot
\sum_{\X_+} |r(x)-s(x)|^2
\leq
\frac 1 {\a_s}
\bigl(\sum_{\X} |r(x)-s(x)|\bigr)^2
=
\frac 4 {\a_s}\cdot |r-s|^2.
$$
\enddemo

We proceed to the proof of Theorem 1.
Let $\pi=\Lw(Y^n, X^n)$ be a coupling of $p^n=\Lw(Y^n)$ and  $q^n=\Lw(X^n)$ that
achieves $W_2(p^n,q^n)$.
\bsk

We apply induction on $n$. Assume that the theorem holds for $n-1$.
\bsk

By the expansion formula for relative entropy we have
$$
D(p^n||q^n)=\frac 1 n\cdot \sumi
D(Y_i||X_i)+\frac 1 n\cdot\sumi D(\bar Y_i|Y_i||\bar q_i|Y_i).
\tag 2.2
$$
For each fixed $y_i$, the measure $\bar q_i(\cdot|y_i)$ satisfies the conditions of the theorem.
By the induction hypothesis,
$$
\align
&\frac 1 n\cdot \sumi D(\bar Y_i|Y_i||\bar q_i(\cdot|Y_i))\\
&\leq
\frac 1 n\cdot \frac {4C} \a\cdot
\sumi \sum_{j\neq i}\bigl|p_j(\cdot|\bar Y_j)-Q_j(\cdot|\bar Y_j)\bigr|^2\\
&=
\biggl(1-\frac 1 n\biggr)\cdot\frac {4C} \a
\cdot \sum_{j=1}^n \bigl|p_j(\cdot|\bar Y_j)-Q_j(\cdot|\bar Y_j)\bigr|^2.
\tag 2.3
\endalign
$$
\bsk

To estimate the first term in the right-hand-side of (2.2) ,
observe that the definition of $\a$ implies
that for any $i\in[1,n]$ and $x\in\X$, $Pr\{X_i=x\}\geq \a$.
Thus by Lemma 2 we have
$$
D(Y_i||X_i)
\leq
\frac 4 \a\cdot \bigl|\Lw(Y_i)-\Lw(X_i)\bigr|^2.
\tag 2.4
$$
Further, condition (1.3) implies
$$
\align
&\sumi
\bigl|\Lw(Y_i)-\Lw(X_i)\bigr|^2
\leq
\sumi
Pr_\pi^2\{Y_i\neq X_i\}=W_2^2(p^n,q^n)\\
&\leq
C
\cdot
\Bbb E
\sumi
\bigl|p_i(\cdot|\bar Y_i)-q_i(\cdot|\bar Y_i)\bigr|^2.
\tag 2.5
\endalign
$$
Putting together (2.4) and (2.5),
it follows that the first term on the right-hand-side of (2.2)
can be bounded as follows:
$$
\frac 1 n\cdot \sumi D(Y_i||X_i)
\leq
\frac 1 n\cdot \frac {2C} \a
\cdot
 \sumi
\Bbb E \bigl|p_i(\cdot|\bar Y_i)-q_i(\cdot|\bar Y_i)\bigr|^2.
\tag 2.6
$$
Substituting (2.3) and (2.6) into (2.2) we get the first \ineq\ in (1.4).
The second \ineq\ follows from the Pinsker-Csisz\'ar-Kullback \ineq.
$\qquad\qquad\qquad\qed$
\bsk
\bsk

\beginsection 3. Proof of Propostion 3.

Let  both $p^n$ and $q^n$ be fixed. We want to show that
(1.3) holds with  $C= 1/ {\bigl(1-||A||\bigr)^2}$, where $A$ is the coupling matrix
for $q^n$. It is enough to prove this for $I=[1,n]$,
since for any $I\su [1,n]$ and $\bar y_I$
the \cd\ $q_I(\cdot|\bar y_I)$ \sats\ Dobrushin's uniqueness condition with a minor of $A$ as
a coupling matrix.
(The idea of the proof for $I=[1,n]$  goes back to Dobrushin's papers [D1], [D2],
 although he worked with another matrix norm.)
\bsk

We are going to  prove that Dobrushin's uniqueness \cond\
implies that  the
Gibbs sampler $\Ga$ is a contraction \wre\ the $W_2$-distance
with rate $1-1/n\cdot (1-||A||)$.
\bsk

To achieve this,
let $r^n$ and $s^n$ be two \pms s on $\X^n$, and
let $U^n$ and $Z^n$ be random \sqe s
representing  $r^n$ and $s^n$, respectively. (I.e., $r^n=\Lw(U^n)$, $s^n=\Lw(Z^n)$.)
\bsk

Select   an index $i\in[1,n]$ at random, and define
$$
{U_k}'=U_k, \qquad {Z_k}'=Z_k
\qute {for}\quad k\neq i.
$$
Then define
$\Lw({U_i}',{Z_i}'| \bar U_i=\bar u_i, \bar Z_i=\bar z_i\}$
as that coupling of $q_i(\cdot|\bar u_i)$
and $q_i(\cdot|\bar z_i)$
that  achieves $|q_i(\cdot|\bar u_i)-q_i(\cdot|\bar z_i)|$.
It is clear that $\Lw({U'}^n)=r^n\Ga$, and $\Lw({Z'}^n)=s^n\Ga$.

\bsk

By the definition of the coupling matrix we have
$$
\Pr\{{U_i}'\neq {Z_i}'\}
\leq
(1-1/n)\cdot Pr\{U_i\neq Z_i\}
+ 1/n\cdot
\sum_{k\neq i}
a_{k,i}\cdot Pr\{U_k\neq Z_k\}.
$$
It follows that
$$
\sqrt{\sumi Pr^2\{{U_i}'\neq {Z_i}'\}}
\leq ||B||\cdot
\sqrt{\sumi   Pr^2\{{U_i}\neq {Z_i}\}},
$$
where
$$
B=(1-1/n)\cdot \Id_n+1/n\cdot A.
$$
Thus
$$
\sqrt{\sumi Pr^2\{{U_i}'\neq {Z_i}'\}}
\leq
\biggl(1-\frac 1 n\cdot(1-||A||\biggr)
\cdot
\sqrt{\sumi   Pr^2\{{U_i}\neq {Z_i}\}}.
$$
This proves the contractivity of $\Ga$
with rate $1-1/n\cdot (1-||A||_2)$.
\bsk

By the triangle \ineq\
$$
W_2(p^n,q^n)
\leq
W_2(p^n, p^n\Ga)+W_2( p^n\Ga,q^n).
$$
By  contractivity of $\Ga$, and since $q^n$ is invariant
\wre\ $\Ga$, it follows that
$$
W_2(p^n,q^n)
\leq
W_2(p^n, p^n\Ga)+ \bigl(1-1/n\cdot (1-||A||)\bigr)\cdot W_2(p^n,q^n),
$$
i.e.,
$$
W_2(p^n,q^n)
\leq
\frac n {1-||A||}\cdot W_2(p^n, p^n\Ga).
$$
But it is easy to see that
$$
W_2(p^n, p^n\Ga)
=
\frac 1 n\cdot
\sqrt{\Bbb E\sumi |p_i(\cdot|\bar Y_i)-q_i(\cdot|\bar Y_i)|^2}.
$$
By the last two \ineqs,  (1.3)
(for $I=[1,n]$), and hence Theorem 2, is proved.
$\qquad\qquad\qquad\qquad\qquad\qquad\qquad\qquad\qquad
\qquad\qquad\qquad\qquad\qquad\qquad
\qed$

\beginsection 4. Gibbs \ms s with the strong mixing property.

\beginsection 4.1. Definitions, notation and statement of Theorem 3.

In this section we work with \ms s on $\X^\La$, where $\La$ is a
subset of the $d$-dimensional cubic lattice  $\Zd$.  Most of the time
$\La$ shall be finite.
\bsk

The lattice points in $\Zd$ shall be called sites.
The distance
$\rho$ on  $\Zd$ is
 $$
 \rho(k,i)=\max_\nu|k_\nu-k_\nu)|,\qute{where}\quad
k=(k_1,k_2,\dots,k_d),\quad  i=(i_1,i_2,\dots,i_d).
$$
\bsk

The notation $\La\subset \subset \Zd$   expresses that
$\La$ is a finite subset of $\Zd$.
\bsk

The elements of $\X$ are called spins, and the
 elements of
the set $\X^\La$ ($\La\subset \Zd$, possibly infinite) are called
spin configurations, or just configurations, over  $\La$.
\bsk

We consider an ensemble of \cd s $q_\La(\cdot|x_{\bar\La})$,
where $\La\subset \subset \Zd$, and $\bar\La$  is the complement of $\La$.
We prefer to write  $\bar x_\La$ in place of $x_{\bar \La}$, and, accordingly,
$q_\La(\cdot|\bar x_{\La})$ in place of  $q_\La(\cdot|x_{\bar\La})$.
The \ms\ $q_\La(\cdot|\bar x_{\La})$ is   considered
as the \cd\  of a random
spin configurations  over  $\La$, given the
spin configuration
outside of $\La$.
For a site $i\in\Zd$  we use the notation $q_i(\cdot|\bar x_{ i})$.

\bsk

The \cd\   $q_\La(\cdot|\bar x_{\La})$  ($\La\subset \subset \Zd$,
 $\bar x_{\La}\in \X^{\bar\La}$) naturally defines the
 \cd s $q_M(\cdot|\bar x_M)$ for any $M\subset \La$.
We assume that
the  \cd s $q_\La(\cdot|\bar x_{\La})$ satisfy the natural compatibility
\cond s. The \cd\   $q_\La(\cdot|\bar x_{\La})$ also defines, for $M\su\La$,
the \cd\  $q_M(\cdot|\bar x_{\La})$.
\bsk

If  the compatibility \cond s hold then
there exists at least one probability \ms\ $q=\Lw(X)$
on the space of \config s $\X^{\Zd}$,
compatible with the
 \cd s $q_\La(\cdot|\bar x_\La)$:
$$
\Lw\bigl(X_\La|\bar X_{\La}=\bar x_{\La}\bigr)
=
q_\La(\cdot|\bar x_\La).
$$
Here $X_\La$ denotes the marginal of the random \config\
$X$ for the sites in $\La$,
and $\bar x_\La$   is called an outside \config\ for $\La$.
The \cd s $q_\La(\cdot|\bar x_\La)$ are called the \lsp s
of $q$, and $q$  is called a Gibbs \ms\ compatible with the \lsp s
$q_\La(\cdot|\bar x_\La)$.
\bsk

We say that the ensemble of \lsp s $q_\La(\cdot|\bar x_\La)$
has finite range of interaction $R$ (or is Markov of order $R$)
if $q_\La(\cdot|\bar x_\La)$
only depends on those \coord s $x_k$ ($k\in \bar\La$)
that are in the $R$-neighborhood of $\La$.
\bsk

In general, the \lsp s do not uniquely determine the Gibbs \ms.
The question of uniqueness has been extensively studied
in the case of \lsp s with finite range of interaction,
 and a sufficient condition for uniqueness was given
 by R. Dobrushin and S. Shlosman
[D-Sh1]. But the general question of uniqueness
is open, even for \ms s with finite range of interaction.
\bsk

A property stronger than  uniqueness  is strong mixing.
\bsk

In their celebrated paper [D-Sh2] in 1987, R. Dobrushin and S. Shlosman
gave a characterization of  complete analyticity of Markov
Gibbs \ms s over   $\Zd$. Their characterization
was formulated in twelve
conditions which were proved to be equivalent, and are referred to as
\DSM. The following definition is the same as one of these twelve
(III C),
except that we do not assume Markovity, and replace
the function $K\cdot \exp(-\ga r)$ by a more general function
$\f(r)$. In the Markov case $\f(r)$ necessarily  shall have
the form $K\cdot \exp(-\ga r)$.
\bsk

In order to define strong mixing, let
$\f: \Bbb Z_+\mapsto \Bbb R_+$ be a function  \sat\
$$
\sum_{i\in \Zd}\varphi\bigl(\rho(0,i)\bigr)<\infty.
\tag 4.1.1
$$

\definition{Definition:  Strong mixing}
The Gibbs \ms\  $q$ is called strongly mixing
with coupling function $\f$
if  for any sets $M\su\La\susu\Zd$
and any two outside \config s
$\bar y_\La$ and $\bar z_\La$
differing only at one single site $k\notin \La$:
$$
\bigl|
q_M(\cdot|\bar    y_\La)
-
q_M(\cdot|\bar    y_\La)
\bigr|
\leq
\f(\rho(k,M)\bigr).
\tag 4.1.2
$$
\enddefinition
\bsk

For stationary \pms s on $\X^{\Bbb Z}$, this definition is more restrictive than usual strong mixing,
and is equivalent to  $\Phi$-mixing.
On $\Zd$ the term strong mixing has been only used
for Markov fields, and for simplicity we extend its use  without adding any qualification.
\bsk

 Our aim in this section is to prove the following
\bsk

\proclaim{Theorem 3}
\nli
Assume that
the  ensemble  $q_\La(\cdot|\bar x_\La)$
\sats\ the \SM\ with coupling function  $\f$.
Moreover, assume that
$$
\a\triangleq\inf q_i(x_i|\bar x_i)>0,
$$
where the infimum is taken for all $x\in\X^{\Zd}$ and $i\in\Zd$
such that $q_i(x_i|\bar x_i)>0$.
Then,  for fixed $\La\susu\Zd$ and outside \config\ $\bar y_\La$,
the \cd\
$q_\La(\cdot|\bar y_\La)$,
as a \ms\  on $\X^\La$, \sats\
condition (1.3) of Theorem 1, with a constant $C$, independent of
$\La$ and $\bar y_\La$.
Moreover, it is enough to assume (4.1.2) for sets $\La$ of
diameter at most $m_0$, where
$m_0$ depends on the dimension $d$ and the function $\f$.
The constant $C$ depends on the dimension $d$, the function $\f$ and on
 $\a$.
\endproclaim
\bsk

\remark{Remark} If $q$ has finite range of interaction then
 Theorem 3 implies that
condition (4.1.2) is
constructive, in the sense of Dobrushin and Shlosman.
\endremark
\bsk

There is another approach to strong mixing, for \ms s $q$ on
$\X^{\Zd}$ with finite range of interaction.
This approach was developed by
 E. Olivieri, P. Picco and F. Martinelli; c.f. [M-O1].
Their aim was to replace the above condition of strong mixing ((4.1.2))
by a milder one, requiring (4.1.2) only for
"non-pathological" sets $\La$, i.e. for sets
whose boundary is much smaller then their volume.
Martinelli and Olivieri [M-O2] proved a \LSI\ under this modified
condition,
for \ms s $q$ with finite range of interaction.
 In Appendix B we briefly
sketch the  Olivieri-Picco-Martinelli approach, and how to modify
Theorem 1 and the Auxiliary Theorem (below), to get
\LSIs\  under this weaker assumption.

\beginsection 4.2. Proof of Theorem 3

Consider the infinite symmetric matrix
$$
\Phi=
\biggl(\f\bigl(\rho(k,i)\bigr)\biggr)_{k,i\in\Zd}.
$$
Since the entries  are non-negative,
and  the row-sums   equal,
$||\Phi||$ equals  the row-sum:
$$
||\Phi||
=
\sum_{i\in \Zd}\f(\rho(0,i)).
$$

\bsk

Fix a $\La\susu\Zd$, an outside \config\ $\byL$
and a
$p_\La\in  \PXLa$.
It is enough to prove that
$$
W_2^2\bigl(  p_\La,q_\La\bigl(\cdot|\byL)\bigr)
\leq
C\cdot
\Bbb E\sum_{i\in \La}
W_2^2
\bigl(
p_i(\cdot|\bar Y_i),
q_i(\cdot|\bar Y_i)
\bigr),
\tag 4.2.1
$$
(with $C$ independent of $\La$ and $\bar y_\La$),
since for any $M\su \La$ and any fixed $ y_{\La\setm M}$,
the \cd\ $q_M(\cdot|\bar y_M)$
(where $\bar y_M= (y_{\La\setm M},\bar y_\La)$)
\sats\ the  \SM\ with the same
function $\f$.
\bsk

We start with a weaker version of (4.2.1).
\bsk

\definition{Notation}
\nli
Let $\I_m=\I_m(\La)$ denote the set of $m$-sided cubes in $\Zd$
that intersect $\La$. Set
$$
\Theta_m
\triangleq
\min_{R}\biggl[||\Phi||\cdot \frac{d\cdot R} m
+
2d\cdot\sum_{r=R}^\infty
(2r+1)^{d-1} \f(r)\biggr].
\tag 4.2.2
$$
\enddefinition

Note that we can achieve
$$
 \Theta_m<1,
\tag 4.2.3
$$
by
selecting  $R$ large enough to make
the  second term in (4.2.2)  small, and then selecting $m$.
\bsk

\proclaim{Auxiliary Theorem}
If $m$ is so large that $ \Theta_m<1$
then
$$
\align
&W_2^2\bigl(p_\La,q_\La(\cdot|\bar y_\La)\bigr)\\
&\leq
\frac 1 {m^d}
\cdot
\frac 1
{(1-\Theta_m)^2}
\cdot
\sum_{I\in\I_m}
\Bbb E
W_2^2
\bigl(
p_{I\cap\La}(\cdot|\bar Y_{I\cap \La}),
q_{I\cap\La}(\cdot|\bar Y_{I\cap \La})
\bigr)\\
&\leq
\frac 1
{(1-\Theta_m)^2}
\cdot
\sum_{I\in\I_m}
\Bbb E
\bigl|
p_{I\cap\La}(\cdot|\bar Y_{ I\cap\La})
-
q_{I\cap\La}(\cdot|\bar Y_{ I\cap\La})
\bigr|^2.
\tag 4.2.4
\endalign
$$
If the ensemble $q_\La(\cdot|\bar x_\la)$ has finite range of interaction $R$
then
 the Auxiliary Theorem holds with
$||\Phi||\cdot\frac{d\cdot R} m $ in place of $\Theta_m$.
\endproclaim
\bsk

The second \ineq\ in (4.2.4) follows from the first one by   the trivial \ineq\
$$
W_2^2(r^n,s^n)\leq n\cdot \bigl|r^n-s^n\bigr|^2\qute{for}\quad
r^n,s^n\in\PXn.
$$
\bsk

The proof of the  Auxiliary Theorem follows  that
of Theorem 2, but
we use a more general  Gibbs sampler,
updating (the intersection of $\La$ with)
an $m$-sided cube at a time, not just one site.
Let us extend the definition of $p_\La$ so that on
$\bar \La$
it be concentrated on the fixed $\byL$.
\bsk

\definition{Definition}
\nli
For $I\in\I_m$ let
$\Ga_I: \PP(\X^\La)\mapsto \PP(\X^\La)$ be the  \Mk:
$$
\Ga_I(z_\La|y_\La)
=
\de_{y_{\La\setm I},z_{\La\setm I}}
\cdot
q_{I\cap\La}(z_{I\cap\La}|\bar y_{I\cap \La}).
$$
(For $k\in \bar \La$,
 $y_k$ is defined by the fixed  $\byL$).
Then set
$$
\Ga_{\I_m}
=\frac 1 {|\I_m|} \cdot \sum_{I\in\I_m} \Ga_I.
$$
Then $\Ga_{\I_m}$ preserves, and is reversible
\wre, $q_\La(\cdot| \byL)$.
We call $\Ga_{\I_m}$ the Gibbs sampler for  \ms\ $q_\La(\cdot|\byL)$,
defined by
the \lsp s
$q_{I\cap\La}(\cdot| \bar y_{I\cap \La})$,
$I\in \I_m$.
\enddefinition

\bsk

\demo{Proof of the Auxiliary Theorem}
\nli
To estimate $W_2^2\bigl(p_\La,q_\La(\cdot| \byL)\bigr)$, we are going to  prove
that if (4.2.3) holds then the
Gibbs sampler $\Ga_{\I_m}$ is a contraction \wre\ the $W_2$-distance,
with rate $1-m^d/|\I_m|\cdot (1-\Theta_m)$.
\bsk

To achieve this,
let $r$ and $s$ be two \pms s on $\X^\La$, and
let $Y$ and $Z$ be \rv s
representing $r$ and $s$, respectively. (I.e., $r=\Lw(Y)$, $s=\Lw(Z)$.)
Let the coupling
$\Lw(Y,Z)$  of $r$ and $s$
achieve $W_2(r,s)$.
We extend the definition of $\Lw(Y,Z)$,
letting
$\bar Y_\La =\bar Z_\La = \byL$, where  $\byL$
 is the fixed outside \config.
Let $Y'$ and $Z'$ be \rv s representing  $r\Ga_{\I_m}$ and $s\Ga_{\I_m}$.
\bsk

Suppose that,  when carrying out one step in the Gibbs sampler
 $\Ga_{\I_m}$,
 we have selected a certain $I\in\I_m$.
Then  we can assume that
$$
 {Y_i}'=Y_i\qute{and}\quad {Z_i}'=Z_i\qute{for all}\quad i\in\La\setm I.
$$
Moreover,
$$
\Lw\bigl(
{Y_{I\cap\La}}'
\bigm|
Y_{\La\setm I}=y_{\La\setm I}\bigr)
=
q_{I\cap\La}\bigl(\cdot|\bar y_{I\cap\La}\bigr),
$$
and
$$
\Lw\bigl(
{Z_{I\cap\La}}'
\bigm|
Z_{\La\setm I}=z_{\La\setm I}\bigr)
=
q_{I\cap\La}\bigl(\cdot|\bar z_{I\cap\La}\bigr).
$$
\bsk

At this point we need the following
\bsk

\proclaim{Lemma 3} (The proof is in  Appendix A.)
\nli
Let us fix the set $M\susu\Zd$, together with  two
outside \config s
 $\bar y_M $ and  $\bar z_M $,
differing only at site $k\notin  M$.
Let $Y$ and $Z$ be \rv s realizing
$q_M\bigl(\cdot|\bar y_M \bigr)$ and $q_M\bigl(\cdot|\bar z_M \bigr)$.
Define
$$
J_i=J_{k,M,i}
=
\bigl\{j\in M: \rho(k,j)\ge \rho(k,i)\bigr\}.
\tag 4.2.5
$$
Then there exists a coupling
$\pi=\Lw\bigl(Y,Z|\bar y_M,  \bar z_M \bigr)$
 of $\Lw(Y)$ and $\Lw(Z)$,
\sat
$$
Pr_\pi \bigl\{Y_i\neq Z_i\bigr\}
=
\bigl|q_{J_i}(\cdot|\bar y_M)-q_{J_i}(\cdot|\bar z_M)\bigr|,
\quad i\in M.
$$
If $q$ \sats\ the \SM\ with function $\f$ then, for this coupling,
$$
Pr_\pi \bigl\{Y_i\neq Z_i\bigr\}
\leq
\f\bigr( \rho(k,i)\bigr)\qute{for all}\quad i\in M.
$$
\endproclaim
\bsk

By Lemma 3, for fixed $I$, $\bar y_{I\cap\La}$ and
$\bar z_{I\cap\La}$, we can define a coupling
$$
\align
&\pi_{I\cap\La}
\bigl(\cdot|\bar y_{I\cap\La},  \bar z_{I\cap\La}\bigr)\\
&=
\Lw\bigl({Y_{I\cap\La}}', {Z_{I\cap\La}}'
\bigm|
\bar Y_{I\cap\La}=\bar y_{I\cap\La},
\bar Z_{I\cap\La}=\bar z_{I\cap\La} \bigr),
\endalign
$$
\sat\
$$
\align
&Pr_{\pi_{I\cap\La}}
\bigl\{
{Y_i}'\neq {Z_i}'
\bigm|
\bar Y_{I\cap\La}=\bar y_{I\cap\La},
\bar Z_{I\cap\La}=\bar z_{I\cap\La} \bigr)\\
&\leq
\sum_{k\in \La\setm I}
\de(y_k,z_k)\cdot \f\bigl(\rho(k,i)\bigr),
\qute{for all}\quad  i\in I\cap\La.
\endalign
$$
Thus
$$
\align
&Pr\bigl\{{Y_i}'\neq {Z_i}'\bigm| I\qute{selected}\bigr\}\\
&\leq
\sum_{k\in\La\setm I}
Pr\{Y_k\neq Z_k\}\cdot \f\bigl(\rho(k,i)\bigr)
\qute{for all}\quad  i\in I\cap\La.
\tag 4.2.6
\endalign
$$
\bsk

We calculate $Pr\{{Y_i}'\neq {Z_i}'\}$ by
 averaging  for $I\in \I_m$.
Set $N=|\I_m|$. Since each $i\in \La$ is covered by exactly $m^d$ cubes from $\I_m$,
 ( 4.2.6) implies
$$
\align
&Pr\bigl\{{Y_i}'\neq {Z_i}'\bigr\}\\
&\leq
\bigl(1-\frac {m^d} N\bigr)\cdot Pr\{Y_i\neq Z_i\}
+
\frac 1 N
\cdot
\sum_{I\ni i}\sum_{k\in \La\setm I} Pr\{Y_k\neq Z_k\}
\cdot \f\bigl(\rho(k,i)\bigr).
\tag 4.2.7
\endalign
$$
\bsk

Consider  the vectors
$$
u=\biggl(Pr\{Y_k\neq Z_k\}\biggr)_{k\in\La}\qand
v=\biggl(Pr\bigl\{{Y_i}'\neq {Z_i}'\bigr\}\biggr)_{i\in\La},
$$
and let $D$ denote the matrix with entries
$$
d_{k,i}
=
\f\bigl(\rho(k,i)\bigr)
\cdot
\sum_{I\ni i,\La\setm I\ni k} 1,
\quad k,i\in\La.
$$
\bsk

With this notation, (4.2.7) means that
$$
v\leq
\biggl(
\bigl(1-\frac {m^d} N\bigr)\cdot\text{Id}
+\frac 1 N \cdot D\biggr)\cdot u
$$
coordinatewise, thus
$$
||v||\leq
\biggl(
\bigl(1-\frac {m^d} N\bigr)+||D||
\biggr)
\cdot ||u||.
\tag 4.2.8
$$
\bsk

We claim that
$$
\sum_{I: k\notin I, I\ni i} 1
\leq
\min \bigl\{d\cdot m^{d-1}\cdot\rho(k,i), \quad\ m^d\bigr\}.
$$
Indeed, there are $d$ lattice-hyperplanes separating
$k$ and $i$, and there is exactly one among these that
intersects the line segment (in $\Bbb R^d$)
connecting $k$ and $i$. These facts imply that
an $m$-sided cube can be placed in at most
$d\cdot m^{d-1}\cdot\rho(k,i)$
ways so as to satisfy both conditions $ k\notin I$ and $ I\ni i$.
It follows that
$$
d_{k,i}
\leq
m^d \cdot\f\bigl(\rho(k,i)\bigr)
\cdot
\min \biggl\{\frac{d\cdot\rho(k,i)} m, 1\biggr\}.
\tag 4.2.9
$$
Since the right-hand-side of (4.2.9) is symmetric in $k$ and $i$,
we have
$$
||D||\leq
m^d\cdot\sum_i
\f\bigl(\rho(k,i)\bigr)
\cdot
\min \biggl\{\frac{d\cdot\rho(k,i)} m, 1\biggr\}.
\tag 4.2.10
$$
Now fix an $R$, and
divide the sum  in (4.2.10) into two parts, for $i$ \sat\
$\rho(k,i)\leq R$ and $(\rho(k,i)> R$, respectively.
We see that
$$
||D||
\leq m^d\cdot
\biggl(
\frac{||\Phi||\cdot d\cdot R}m
+
\sum_{\rho(k,i)>R}\f\bigl(\rho(k,i)\bigr)\biggr).
$$
Taking minimum in $R$, we get
$$
||D||\leq m^d\cdot\Theta_m.
\tag 4.2.11
$$

By (4.2.8) and the definition of the vectors $u$ and $v$,
(4.2.11) implies that
$$
\sqrt{\sum_{i\in\La}
Pr^2\bigl\{{Y_i}'\neq {Z_i}'\bigr\}}
\leq
\biggl(1- \frac {m^d} N\cdot (1-\Theta_m)\biggr)
\cdot
\sqrt{
\sum_{k\in\La}
Pr^2\{Y_k\neq Z_k\}},
$$
i.e.,
$$
W_2\bigl( r\Ga_{\I_m}, s\Ga_{\I_m}\bigr)
\leq
\biggl(1- \frac {m^d} N\cdot (1-\Theta_m)\biggr)
\cdot
W_2(r,s).
\tag 4.2.12
$$
The stated contractivity is proved.
\bsk

By the triangle \ineq\ it follows that
$$
\align
&W_2\bigl(p_\La,q_\La(\cdot|\bar y_\La)\bigr)
\leq
W_2\bigl(p_\La,p_\La\Ga_{\I_m}\bigr)
+
W_2\bigl(p_\La\Ga_{\I_m},q_\La(\cdot| \bar y_\La)\bigr)\\
&\leq
W_2\bigl(p_\La,p_\La\Ga_{\I_m}\bigr)
+
\biggl(1- \frac {m^d} N\cdot (1-\Theta_m)\biggr)
\cdot
W_2\bigl(p_\La, q_\La(\cdot|\bar y_\La)\bigr),
\endalign
$$
whence
$$
W_2\bigl(p_\La,q_\La(\cdot|\bar y_\La)\bigr)
\leq
\frac N {m^d}\cdot
\frac 1 {(1-\Theta_m)}
\cdot W_2\bigl(p_\La,p_\La\Ga_{\I_m}\bigr).
\tag 4.2.13
$$
\bsk

To complete the proof if the Auxiliary Theorem, we have to
estimate $ W_2\bigl(p_\La,p_\La\Ga_{\I_m}\bigr)$
in terms of the quantities
\nli
$$
\Bbb E
W_2^2
\bigl(
p_{I\cap\La}(\cdot|\bar Y_{I\cap\La}),
q_{I\cap\La}(\cdot|\bar Y_{I\cap\La}) \bigr).
$$
To do this, fix an $I\in\I_m$, together with a \sqe\ $y_{\La\setm I}\in \X^{\La\setm I}$, and
define a coupling  $\pi_{I\cap\La}(\cdot| y_{\La\setm I})$ of
$p_{I\cap\La}(\cdot|\bar y_{ I\cap\La}) $       and
$q_{I\cap\La}(\cdot|\bar y_{ I\cap\La}) $
that achieves  $W_2$-distance.
We extend  $\pi_{I\cap\La}(\cdot| y_{\La\setm I})$
to a \ms\ on $\X^\La\times \X^\La$  concentrated
on the diagonal
$(y_{\La\setm I}, y_{\La\setm I})$,
for \coord s
outside of $I$.
Finally, we define the coupling $\pi$ of
$p_\La$ and $p_\La\Ga_{\I_m}$ by averaging the distributions
$\pi_{I\cap\La}(\cdot| y_{\La\setm I})$ \wre\
 $I$ and  $y_{\La\setm I}$.

\bsk

Using this construction, an easy computation (using the Cauchy-Schwarz
\ineq)
shows that
$$
W_2^2\bigl(p_\La,p_\La\Ga_{\I_m}\bigr)
\leq
\frac{m^d}{N^2}
\sum_{I\in\I_m}
\Bbb E
W_2^2
\bigl(
p_{I\cap\La}(\cdot|\bar Y_{I\cap\La}),
q_{I\cap\La}(\cdot|\bar Y_{I\cap\La})
\bigr).
\tag 4.2.14
$$
\bsk

Substituting (4.2.14) into (4.2.13), we get the first \ineq\ in (4.2.4). Understanding the proof one easily sees
that the statement
for Gibbs \ms s with finite range of interaction holds true.
The Auxiliary Theorem is proved.
$\qquad\qquad\qquad\qed$
\enddemo
\bsk

To complete the proof of Theorem 3
we have to deduce (4.2.1) from the Auxiliary Theorem.
To do this we need the following
\bsk

\proclaim{Lemma 4} (The proof is in  Appendix A.)
\nli
Let $p^n=\Lw(Y^n)$ and $q^n$ be two \ms s on $\X^n$.
Let $\a$ be defined by (1.1).
Then
$$
\bigl|p^n- q^n\bigr|^2
\leq
\Biggl(\frac 2 {(|\X|\cdot \a)^2}\Biggr)^{n+\log_2 n}
\cdot
\sumi
\Bbb E
\bigl|p_i(\cdot|\bar Y_i)-q_i(\cdot|\bar Y_i)\bigr|^2.
$$
\endproclaim
\bsk
\bsk

Using Lemma 4, we estimate  the terms in the last sum in (4.2.5). We get
$$
\align
&W_2^2\bigl(p_\La,q_\La(\cdot|\bar y_\La)\bigr)\\
&\leq
\frac {m^d}
{(1-\Theta_m)^2}
\cdot
\Biggl(\frac 2 {(|\X|\cdot \a)^2}\Biggr)^{m+\log_2 m}
\cdot
\Bbb E
\sum_{i\in \La}
\bigl|p_i(\cdot| Y_{\La\setm i})-q_i(\cdot|  Y_{\La\setm i},
\bar y_\La)\bigr|^2.
\endalign
$$
Thus (4.2.1)
is fulfilled with
$$
C=\frac {m^d}
{(1-\Theta_m)^2}
\cdot
\Biggl(\frac 2 {(|\X|\cdot \a)^2}\Biggr)^{m+\log_2 m},
$$
as soon as $m$ is large enough for $\Theta_m<1$.
\bsk

We used the \SM\ (4.1.2) in proving Lemma 3, and Lemma 3
was used for subsets of  $m$-sided cubes. It was enough
to consider $m$-sided cubes with
 $m$ so large that $\Theta_m<1$ holds,
a condition depending on $d$ and $\f$. This proves
the last two statements of Theorem 3. $\qquad\qquad\qquad\qquad\qquad\qquad\qquad\qquad\qquad\quad\qed$
\bsk

\remark{Remark} An argument similar to the use of Lemma 4 was also there
in [S-Z2].
\endremark
\bsk

\beginsection Acknowledgement

The author thanks M. Raginsky for providing a simple proof
of Lemma 1.
\bsk

\beginsection Appendix A

\demo{Proof of Lemma 1} (This proof was suggested to the author by M. Raginsky [R].)
\bsk

We use the notions of Hellinger \dis\ and Hellinger affinity:
$$
H(r,s)
=
\biggl(
\sum_{x\in\X}\biggl|\sqrt{r(x)}-\sqrt{s(x)}\biggr|^2
\biggr)^{1/2}
\qand
A(r,s)
=
\sum_{x\in\X}\sqrt{r(x)\cdot s(x)}.
$$
\bsk

The statement of the lemma can be formulated as
$$
|r-s|^2\leq 1-A^2(r,s).
\tag A.1
$$
\bsk

To prove (A.1),
we use the identity
$$
H^2(r,s)
=
2\bigl(1-A(r,s)\bigr).
$$
(A.1) is now
proved by the following chain of equalities and \ineqs:
$$
\align
&|r-s|^2
=
\biggl(
\frac 1 2\cdot\sum_{x\in\X}\bigl|r(x)-s(x)\bigr|
\biggr)^2\\
&=
\frac 1 4
\biggl(
\sum_{x\in\X}
\biggl|\sqrt{r(x)}-\sqrt{s(x)}\biggr|\cdot
\biggl|\sqrt{r(x)}+\sqrt{s(x)}\biggr|
\biggr)^2\\
&\leq
\frac 1 4\cdot
\sum_{x\in\X}
\biggl|\sqrt{r(x)}-\sqrt{s(x)}\biggr|^2
\cdot \sum_{x\in\X}
 \biggl|\sqrt{r(x)}+\sqrt{s(x)}\biggr|^2\\
&=
H^2(r,s)\cdot 2\bigl(1+A(r,s)\bigr)\\
&=
\bigl(1-A(r,s)\bigr)\cdot \bigl(1+A(r,s)\bigr)=1-A^2(r,s).
\endalign
$$
(The inequality follows from the Cauchy-Schwarz \ineq.)
$\qquad\qquad\qquad\qquad\qquad\qquad\qquad\qed$
\enddemo

\demo{Proof of Lemma 3}
\nli
Order the elements of $\La$ so that
$$
\rho(k,i_1)\leq \rho(k,i_2)\leq\dots\leq \rho(k,i_{|\La|}),
$$
i.e., the \sqe\ of sets $J_i=J_{k,M,i}$ (c.f. (4.1.2))
is decreasing in $i$.
Let $Y_{J_i}$ and
$Z_{J_i}$ denote the marginals of $Y$ and $Z$,
respectively, for the sites in $J_i$.
Then
$(Y_{J_1},Y_{J_2},\dots, Y_{J_{|M|}})$ and
$(Z_{J_1},Z_{J_2},\dots, Z_{J_{|M|}})$
are Markov chains. (In fact, $Y_{J_{i+1}}$ is a function of $Y_{J_i}$.)
Therefore, by a theorem of Goldstein [Go],
there exists a coupling
$\pi=\Lw\bigl(Y,Z|\bar y_M,  \bar z_M\bigr)$
 of $\Lw(Y)$ and $\Lw(Z)$,
\sat\
$$
Pr_\pi \bigl\{Y_{J_i}\neq Z_{J_i}\bigr\}
=
\bigl|\Lw(Y_{J_i})-\Lw(Z_{J_i})\bigr|
=
\bigl|q_{J_i}(\cdot|\bar y_M)-q_{J_i}(\cdot|\bar z_M)\bigr|.
$$
Since $i\in J_i$, and
$\rho(k,i)=\rho(k,J_i) $, the statement of Lemma 4 follows.
$\qquad\qquad\qquad\qed$
\enddemo
\bsk

\demo{Proof of Lemma 4}
\nli
Note first that if $r$ and $s$ are \pms s on $\X$, and $r(x), s(x)\ge \a$
then
$$
|r-s|\leq 1-|\X|\cdot \a.
$$
\bsk

Now consider \ms s $p^2=\Lw(Y_1,Z_2)$ and $q^2$ on a product space $\Cal Y\times \Cal Z$,
where $q_2(z_2|y_1)\ge \a_2$,
and  $q_1(y_1|z_2)\ge \a_1$ for all $y_1,z_2\in  \Cal Y\times \Cal Z$.
Then
$$
\bigl|q_2(\cdot|y_1)-q_2(\cdot|{y_1}')\bigr|
\leq
1-|\Cal Z|\cdot \a_2,\qand
\bigl|q_1(\cdot|z_2)-q_2(\cdot|{z_2}')\bigr|
\leq
1-|\Cal Y|\cdot \a_1
$$
for all $y_1, {y_1}'\in \Cal Y$  and $z_2, {z_2}'\in \Cal Z$.
\bsk

Thus in this case
Dobrushin's uniqueness condition is satisfied
with  a $2\times 2$ coupling matrix,
with entries
$1-|\Cal Y|\cdot \a_1$ and $1-|\Cal Z|\cdot \a_2$ outside the diagonal.
(It does not matter that  $\Cal Y$ and $\Cal Z$ may be different.)
The coupling matrix  has norm
$$
\leq
\max
\bigl\{1-|\Cal Y|\cdot \a_1,  1-|\Cal Z|\cdot \a_2\bigr\}.
$$
\bsk

By the argument proving Theorem 2,
it follows that
$$
\align
&W_2(p^2,q^2)\\
&\leq
\max
\biggl\{\frac 1 {\bigl(|\Cal Y|\cdot \a_1\bigr)^2},
       \frac 1 {\bigl(|\Cal Z|\cdot \a_2\bigr)^2}
\biggr\}
\cdot
\Bbb E
\biggl(
\bigl|p_1(\cdot|Y_2)-q_1(\cdot|Y_2)\bigr|^2
+
\bigl|p_2(\cdot|Y_1)-q_2(\cdot|Y_1)\bigr|^2
\biggr),
\endalign
$$
and, consequently,
$$
\align
&\bigl|p^2-q^2|^2\\
&\leq
\max
\biggl\{\frac 2 {\bigl(|\Cal Y|\cdot \a_1\bigr)^2},
       \frac 2 {\bigl(|\Cal Z|\cdot \a_2\bigr)^2}
\biggr\}
\cdot
\Bbb E
\biggl(
\bigl|p_1(\cdot|Y_2)-q_1(\cdot|Y_2)\bigr|^2
+
\bigl|p_2(\cdot|Y_1)-q_2(\cdot|Y_1)\bigr|^2
\biggr).
\tag A1
\endalign
$$
\bsk

Lemma 4 follows from (A1) by a recursive argument, dividing the index set
into two possibly equal parts of size $\lceil\frac n 2\rceil$ and
$\lfloor\frac n 2\rfloor$, and applying (A1) for the two parts.
Then
$$
\max
\biggl\{\frac 2 {\bigl(|\Cal Y|\cdot \a_1\bigr)^2},
       \frac 2 {\bigl(|\Cal Z|\cdot \a_2\bigr)^2}
\biggr\}
$$
shall be replaced by
$$
\biggl(
\frac 2
{(|\X|\cdot\a)^2}
\biggr)
^{\lceil\frac n 2\rceil}.
$$
Repeating this step about $\log_2 n$ times we get the statement of the lemma.
$\qquad\qquad\qed$
\enddemo
\bsk

\beginsection Appendix B

Let $\Zd/l$ ($l\ge 1$ integer) denote the
sub-lattice in $\Zd$, consisting of points whose  coordinates
are all multiples of $l$, and let $\Cal C_l$ denote the set of 
finite unions of $l$-sided cubes
with vertices in $\Zd/l$.

\bsk


The approach by   Olivieri and Picco
is based on the following definition of strong mixing:
\bsk

\definition{Definition by Olivieri and Picco}
\nli
The Gibbs \ms\  $q$  on $\X^{Z^d}$ with finite range of interaction
is called strongly mixing over $\Cal C_l$,
if there exist numbers  $\ga>0$, $K>0$   such that:
 for any sets $\La\in\Cal C_l$, $M\su\La$
and any two outside \config s
$\bar y_\La$ and  $\bar z_\La$  differing only at a single
site $k\notin \La$,
we have
$$
\bigl|
q_M(\cdot|\bar    y_\La)
-
q_M(\cdot|\bar    z_\La)
\bigr|
\leq
K\cdot \exp\bigl(-\ga\cdot\rho(k,M)\bigr).
\tag B.1
$$
\enddefinition

In force of the following theorem, if $l$ is sufficiently large  then it is enough
to require (B.1) just for cubes in $\Cal C_l$,
to get (B.1) for all $\La\in\Cal C_l$, however,
with a different
$\ga$ and $K$.

\proclaim{Olivieri and Picco's Effectivity Theorem, [O-P], [M-O1]}
\nli
Assume that the Gibbs \ms\ $q$ on $\X^{Z^d}$ has finite range of interaction.
For any $\ga,K>0$ there exists an $l_0$ such that:
if for some $l\ge l_0$ (B.1) holds
for all $l$-sided \it{cubes} $\La\in\Cal C_{l}$, all $M\su\La$ and all
$k\notin \La$, then (B.1) also holds  for all
$\La\in \Cal C_l$,  and $M$ and $k$ as above, with  different
$\ga$ and $K$.
\endproclaim

We use a slightly more general definition, although we cannot justify it
with some analog of the above Effectivity Theorem:
\bsk

\definition{Definition:  Strong mixing over $\Cal C_l$}
\nli
Let
$\f: \Bbb Z_+\mapsto \Bbb R_+$ be a function  \sat\ (4.1.1).
Fix an integer $l\ge 1$.
The ensemble of \cd s $q_\La(\cdot|\bar x_\la)$
 on $\X^{Z^d}$
is called strongly mixing  over $\Cal C_l$,
with coupling function $\f$, if
for any sets  $\La\in\Cal C_l$, $M\su\La$, and any
two outside \config s
$\bar y_\La$ and  $\bar z_\La$  differing only at the single
site $k$,  (4.1.2) holds.
(We do not assume finite range of interaction.)
\enddefinition

For \ms s  strongly mixing over  $\Cal C_l$ one can prove a  \LSI\
by means of the following  modifications of Theorem 1 and the Auxiliary Theorem:
\bsk

\proclaim{Theorem 1'}
\nli
Consider a \ms\ $q^\La$ on $\X^\La=\prod_{j=1}^n \X^{\La_j}$, where
$$
\La=\cup_{j=1}^n \La_j,\quad \La_j\cap \La_k=\emptyset
\qute{for}\quad j\neq k,
\quad |\La_j|=m.
$$
Set
$$
\a=\min \bigl\{q_i(x_i|\bar x_i):\quad q_\La(x_\La)>0,i\in\La\bigr\}.
$$
Fix a $p_\La=\Lw(Y_\La)$ on $\X^\La$ \sat
$$
q_\La(x_\La)=0 \implies p_\La(x_\La)=0.
$$
Assume that
 $q_\La$
\sats\ all
the \ineqs\
$$
W_2^2\bigl(  p_I\bigl(\cdot|\bar y_I\bigr),
q_I\bigl(\cdot|\bar y_I\bigr)\bigr)
\leq
C\cdot
\Bbb E
\biggl\{
\sum_{\La_j\su I}
W^2
\bigl(
 p_{\La_j}(\cdot|\bar Y_{\La_j}),q_{\La_j}(\cdot|\bar Y_{\La_j})\bigr)
\biggm|\bar Y_I=\bar y_I
\biggr\},
$$
where   $I\su\La$  is the
 union of some of the sets $\La_j$,
 and  $\bar y_I\in \X^{\La\setm I}$ is a fixed \sqe.
Then
$$
D(p_\La||q_\La)
\leq
\frac {4Cm} {\a^m}
\cdot
\sum_{j=1}^n
\Bbb E
W^2
\bigl(
 p_{\La_j}(\cdot|\bar Y_{\La_j}),q_{\La_j}(\cdot|\bar Y_{\La_j})
 \bigr).
$$
\endproclaim
\bsk

This can be proved by the same argument as Theorem 1, using
Lemma 1, the \ineqs
$$
\biggl| p_{\La_j}(\cdot|\bar Y_{\La_j})
-q_{\La_j}(\cdot|\bar Y_{\La_j})
\biggr|^2
\leq
m\cdot
W^2
\bigl(
 p_{\La_j}(\cdot|\bar Y_{\La_j}),q_{\La_j}(\cdot|\bar Y_{\La_j})
 \bigr),
$$
and, in each induction step,
 fixing a whole new block
$Y_{\La_j}$. 
\bsk

\proclaim{Auxiliary Theorem for \ms s strongly mixing over $\Cal C_l$}
\nli
Fix an integer $l$, and
assume that the ensemble of \cd s $q_\La(\cdot|\bar x_\La)$
 on $\X^{Z^d}$ \sats\
the \SM\ over  $\Cal C_l$, with coupling function $\f$.
Let $\La\in\Cal C_l$, and fix an outside \config\ $\bar y_\La$.
For fixed $m$
let $\Cal I_{ml}$ denote the set of $m\cdot l$-sided cubes from $\Cal C_l$
intersecting $\La$. Then for large enough $m$ and any \ms\ $p_\La$ on $\X^\La$
$$
\align
&W_2^2\bigl(p_\La,q_\La(\cdot|\bar y_\La)\bigr)\\
&\leq
C
\cdot
\sum_{I\in\I_{ml}}
\Bbb E
W_2^2
\bigl(
p_{I\cap\La}(\cdot|\bar Y_{I\cap \La}),
q_{I\cap\La}(\cdot|\bar Y_{I\cap \La})
\bigr)\\
&\leq
C\cdot m^d\cdot
\sum_{I\in\I_{ml}}
\Bbb E
\bigl|
p_{I\cap\La}(\cdot|\bar Y_{ I\cap\La})
-
q_{I\cap\La}(\cdot|\bar Y_{ I\cap\La})
\bigr|^2,
\endalign
$$
where $C$ and $m$ depend on the dimension $d$ and on the function $\f$.
\endproclaim
\bsk

The proof uses   a Gibbs sampler, updating (intersections with
$\La$ of) randomly chosen
cubes of side $m\cdot l$ from $\Cal C_l$. (For an appropriate $m$.)

\enddocument